# Average Number of Lattice Points in a Disk

Sujay Jayakar    Robert S. Strichartz [*]

**Abstract**

The difference between the number of lattice points in a disk of radius $\sqrt{t}/2\pi$ and the area of the disk $t/4\pi$ is equal to the error in the Weyl asymptotic estimate for the eigenvalue counting function of the Laplacian on the standard flat torus. We give a sharp asymptotic expression for the average value of the difference over the interval $0 \leq t \leq R$. We obtain similar results for families of ellipses. We also obtain relations to the eigenvalue counting function for the Klein bottle and projective plane.

## 1 The simplest case

Consider the standard flat torus $[0,1] \times [0,1]$ with boundaries identified. The eigenfunctions of the Laplacian are $e^{2\pi i n \cdot x}$ for $n \in \mathbb{Z}^2$ with eigenvalues $(2\pi)^2 |n|^2$, so the eigenvalue counting function is

$$N(t) = \#\left\{ n \in \mathbb{Z}^d : |n| \leq \sqrt{t}/2\pi \right\}, \qquad (1.1)$$

the number of lattice points inside the disk $B_{\sqrt{t}/2\pi}$ of radius $\sqrt{t}/2\pi$ about the origin. To first approximation $N(t)$ is the area of the disk $t/4\pi$, and this is exactly the Weyl asymptotic law. The problem of estimating the difference

$$D(t) = N(t) - \frac{t}{4\pi} \qquad (1.2)$$

is notoriously difficult (conjectured to be $O(t^{1/4+\epsilon})$ for every $\epsilon > 0$). Here

---

[*]Research supported in part by NSF grant #DMS-0652440

2010 *Mathematics Subject Classification.* Primary 35J05; Primary 42B99.

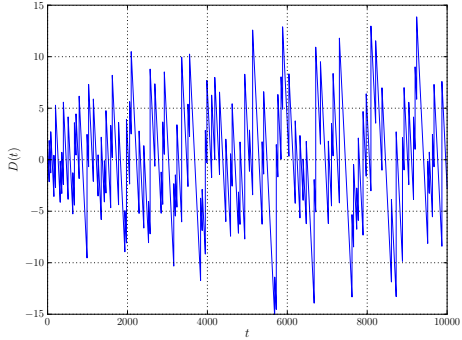
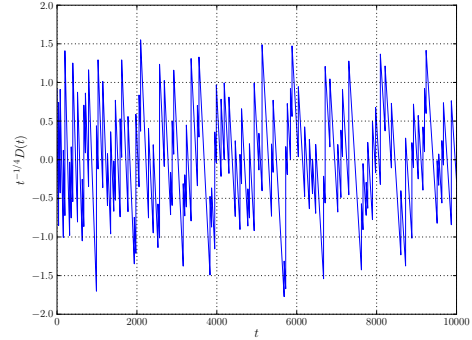

Figure 1: $D(t)$          Figure 2: $t^{1/4}D(t)$

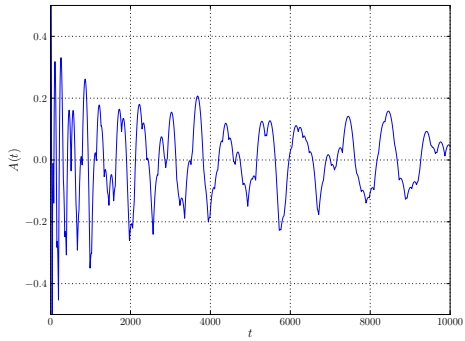
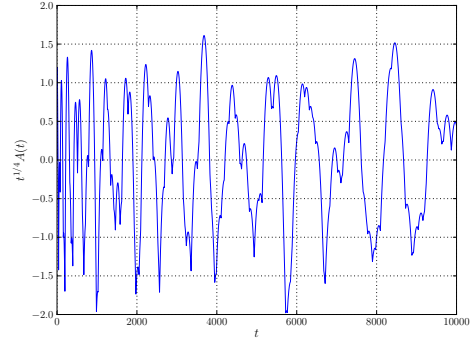

Figure 3: $A(t)$          Figure 4: $t^{1/4}A(t)$

we study the simpler problem of approximating the average value

$$A(R) = \frac{1}{R} \int_0^R D(t)\, dt. \tag{1.3}$$

Note that we are not taking the absolute value of $D(t)$ in the average, so we may exploit the cancellation from regions where $N(t)$ is greater than and less than $t/4\pi$. We will show that $A(R) = O(R^{-1/4})$ as $R \to \infty$, and more precisely

$$A(R) = g(R^{-1/2})R^{-1/4} + O(R^{-3/4}) \quad \text{as } R \to \infty \tag{1.4}$$

where $g(R)$ is an explicit uniformly almost periodic function of mean value zero. Somewhat different but related ideas are given in Bleher [2, 3]. The following Lemma is well-known (see [4], p. 74), but we include the proof for the convenience of the reader.

**Lemma 1.** *We have*

$$A(R) = \sum_{n \neq 0} \frac{1}{\pi |n|^2} J_2(|n|\sqrt{R}) \tag{1.5}$$

*where $n = (n_1, n_2)$ is a variable in $\mathbb{Z}^2$, and $J_2$ denotes the Bessel function. The series in (1.5) converges uniformly and absolutely.*

*Proof.* Let $\chi_t$ denote the characteristic function of the ball $B_{\sqrt{t}/2\pi}$. It is well-known that

$$\hat{\chi}_t(z) = \begin{cases} \frac{\sqrt{t}}{2\pi |z|} J_1(|z|\sqrt{t}) & z \neq 0 \\ \frac{t}{4\pi} & z = 0 \end{cases} \tag{1.6}$$

Following standard methods (see [7] or [4]) we apply the Poisson summation formula to

$$F_{R,\delta} = \frac{1}{R} \int_0^R \chi_t * \psi_\delta \, dt, \tag{1.7}$$

where $\psi_\delta$ is a smooth approximate identity. The $\psi_\delta$ convolution makes $F_{R,\delta}$ smooth, but eventually we will let $\delta \to 0$. Note that

$$\frac{1}{R} \int_0^R N(t) \, dt = \lim_{\delta \to 0} \sum_{n \in \mathbb{Z}^2} F_{R,\delta}(n). \tag{1.8}$$

The Poisson summation formula gives

$$\sum_{n \in \mathbb{Z}^2} F_{R,\delta}(n) = \sum_{n \in \mathbb{Z}^2} \hat{F}_{R,\delta}(n) \tag{1.9}$$

$$= \sum_{n \in \mathbb{Z}^2} \frac{1}{R} \int_0^R \hat{\chi}_t(n) \hat{\psi}(\delta n) \, dt$$

$$= \frac{1}{R} \int_0^R \frac{t}{4\pi} \, dt + \sum_{n \neq 0} \frac{1}{R} \int_0^R \frac{\sqrt{t}}{2\pi |n|} J_1(|n|\sqrt{t}) \, dt \, \hat{\psi}(\delta n)$$

by (1.6). Combining (1.8) and (1.9) yields

$$A(R) = \lim_{\delta \to 0} \sum_{n \neq 0} \frac{1}{R} \int_0^R \frac{\sqrt{t}}{2\pi |n|} J_1(|n|\sqrt{t}) \, dt \, \hat{\psi}(\delta n). \tag{1.10}$$

Now we use the property of Bessel functions ([6])

$$\int_0^R s^{\alpha+1} J_\alpha(s)\, ds = R^{\alpha+1} J_{\alpha+1}(R) \tag{1.11}$$

for $\alpha = 1$, together with the change of variables $s = |n|\sqrt{t}$, to evaluate the integral in (1.10)

$$\frac{1}{R} \int_0^R \frac{\sqrt{t}}{2\pi|n|} J_1(|n|\sqrt{t})\, dt = \frac{1}{R} \int_0^{|n|\sqrt{R}} \frac{s^2}{\pi|n|^4} J_1(s)\, ds \tag{1.12}$$
$$= \frac{1}{\pi|n|^2} J_2(|n|R),$$

and substitute this into (1.10) to obtain

$$A(R) = \lim_{\delta \to 0} \sum_{n \neq 0} \frac{1}{\pi|n|^2} J_2(|n|\sqrt{R}) \hat{\psi}(\delta n). \tag{1.13}$$

The estimate $J_2(|n|\sqrt{R}) = O(\frac{1}{|n|^{1/2} R^{1/4}})$ shows the convergence of the sum in (1.13) without the term $\hat{\psi}(\delta n)$, so we can take the limit in (1.13) and obtain (1.5). □

**Theorem 2.** *Consider the uniformly almost periodic function with mean value zero*

$$g(x) = -\frac{\sqrt{2}}{\pi^{3/2}} \sum_{n \neq 0} |n|^{-5/2} \cos\left(|n|x - \frac{\pi}{4}\right). \tag{1.14}$$

*We have*
$$A(R) = g(R^{1/2}) R^{-1/4} + O(R^{-3/4}) \quad \text{as } R \to \infty. \tag{1.15}$$

*More generally, there exists a sequence of uniformly almost periodic functions $g_1, g_2, \ldots$ with $g_1 = g$ such that for any $n$,*

$$A(R) = \sum_{j=1}^n g_j(R^{1/2}) R^{\frac{1}{4} - \frac{j}{2}} + O(R^{-\frac{1}{4} - \frac{n}{2}}). \tag{1.16}$$

*Proof.* We use the well-known asymptotic expression for Bessel functions

$$J_\alpha(x) = \sqrt{\frac{2}{\pi}} x^{-1/2} \cos\left(x - \frac{1}{2}\alpha\pi - \frac{\pi}{4}\right) + O(x^{-3/2}) \quad \text{as } x \to \infty. \tag{1.17}$$

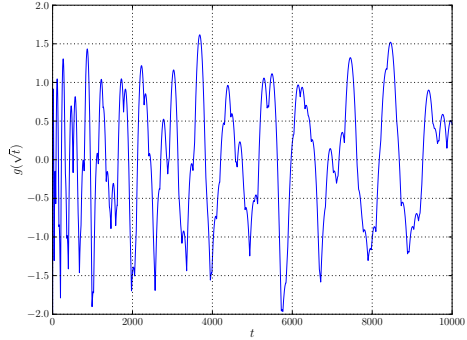 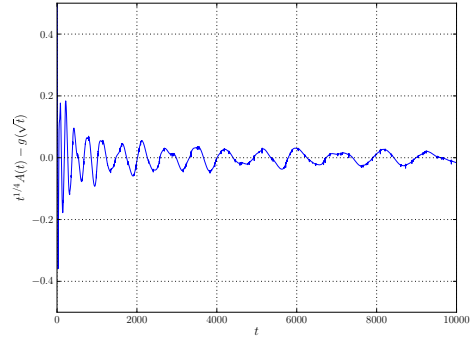

Figure 5: $g(\sqrt{t})$  Figure 6: $t^{1/4}A(t) - g(\sqrt{t})$

When $\alpha = 2$ this is

$$J_2(x) = -\sqrt{\frac{2}{\pi}}x^{-1/2}\cos\left(x - \frac{\pi}{4}\right) + O(x^{-3/2}), \qquad (1.18)$$

and we substitute this into (1.5) with $x = |n|\sqrt{R}$ to obtain

$$A(R) = g(R^{1/2})R^{-1/4} + \sum_{n\neq 0}\frac{1}{|n|^2}O((n\sqrt{R})^{-3/2}). \qquad (1.19)$$

It is easy to see that the remainder term in (1.19) is $O(R^{-3/4})$, so (1.19) yields (1.15). To obtain the more refined asymptotic expression (1.16) we use the known more refined asymptotic expansion for Bessel functions (see [6]). In particular we note that it is possible to obtain explicit series expansions of the functions $g_j$; for example,

$$g_2(x) = \frac{15\sqrt{2}}{8\pi^{3/2}}\sum_{n\neq 0}|n|^{-7/2}\sin\left(|n|x - \frac{\pi}{4}\right). \qquad (1.20)$$

□

It is also reasonable to consider the function $N((2\pi r)^2)$ that counts the number of lattice points inside the ball $B_r$ of radius $r$, the difference $D((2\pi r)^2) = N((2\pi r)^2) - \pi r^2$, and the average with respect to the radius variable

$$\tilde{A}(R) = \frac{1}{R}\int_0^R D((2\pi r)^2)\,dr. \qquad (1.21)$$

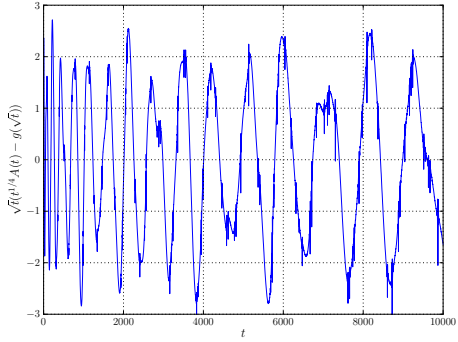
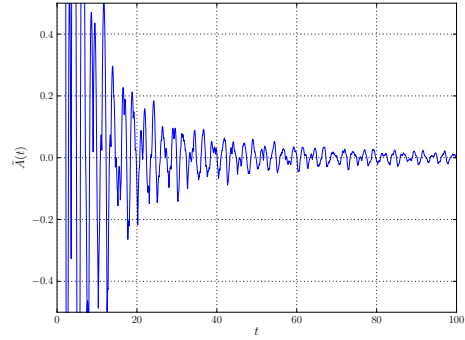

Figure 7: $\sqrt{t}(t^{1/4}A(t) - g(\sqrt{t}))$

Figure 8: $\tilde{A}(t)$

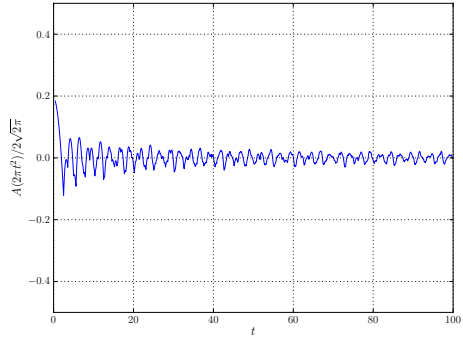
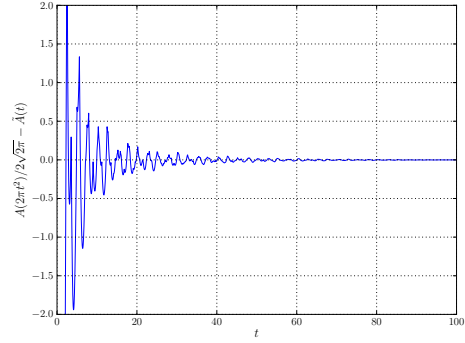

Figure 9: $\frac{1}{2\sqrt{2\pi}}A(2\pi t^2)$

Figure 10: $\frac{1}{2\sqrt{2\pi}}A(2\pi t^2) - \tilde{A}(t)$

This is a different average, but a change of variable shows that

$$\tilde{A}(R) = \frac{1}{2}R\left(\frac{1}{2\pi R^2}\int_0^{2\pi R^2} D(t)\frac{dt}{t^{1/2}}\right). \quad (1.22)$$

Since most of the contribution to the integral occurs for values of $t$ near $2\pi R^2$, we see that $\tilde{A}(R)$ has the same asymptotics as $\frac{1}{2\sqrt{2\pi}}A(2\pi R^2)$.

In Figure 1 we show the graph of $D(t)$ and in Figure 2 the graph of $t^{-1/4}D(t)$. This illustrates the rough $t^{1/4}$ growth rate of $D(t)$. In Figure 3 we show the graph of $A(t)$, and Figure 4 the graph of $t^{1/4}A(t)$. Figure 5 shows the graph of $g(\sqrt{t})$, which is almost identical to Figure 4 for large $t$. Figure 6 shows the difference of $t^{1/4}A(t)$ and $g(\sqrt{t})$, and Figure 7 shows this difference multiplied by $t^{1/2}$. Figure 8 shows the graph of $\tilde{A}(t)$. Figure 9 shows the graph of $\frac{1}{2\sqrt{2\pi}}A(2\pi t^2)$, which agrees with Figure 8 for large $t$, and Figure 10 shows the difference. For more data see the website [5].

## 2 The general case

Consider the general flat 2-dimensional torus, $\mathbb{R}^2/\mathcal{L}$ for some lattice $\mathcal{L}$. The eigenfunctions of the Laplacian (restriction of the standard $\mathbb{R}^2$ Laplacian) have the form $e^{2\pi i x \cdot \xi}$ for $\xi$ in the dual lattice $\mathcal{L}'$, with eigenvalues $(2\pi)^2 \xi \cdot \xi$. By diagonalizing the quadratic form $\xi \cdot \xi$ on $\mathcal{L}'$ we can find an orthonormal basis $v_1, v_2$ in $\mathbb{R}^2$ and positive constants $a_1, a_2$, such that the eigenvalues are

$$(2\pi)^2 \left( \left( \frac{n \cdot v_1}{a_1} \right)^2 + \left( \frac{n \cdot v_2}{a_2} \right)^2 \right) \quad \text{for } n \in \mathbb{Z}^2.$$

Thus the eigenvalue counting function is

$$N(t) = \# \left\{ n \in \mathbb{Z}^2 : \left( \left( \frac{n \cdot v_1}{a_1} \right)^2 + \left( \frac{n \cdot v_2}{a_2} \right)^2 \right)^{1/2} \leq \sqrt{t}/2\pi \right\}. \quad (2.1)$$

In place of disks we consider the family of ellipses

$$E_t = \left\{ x \in \mathbb{R}^2 : (2\pi)^2 \left( \left( \frac{x \cdot v_1}{a_1} \right)^2 + \left( \frac{x \cdot v_2}{a_2} \right)^2 \right) \leq t \right\}. \quad (2.2)$$

Of course $N(t)$ is just the number of lattice points in $E_t$, and the volume of $E_t$ is $\frac{a_1 a_2 t}{4\pi}$. Again we write $D(t) = N(t) - \frac{a_1 a_2 t}{4\pi}$ for the difference and define the average $A(R)$ by (1.3). The analog of (1.6) is

$$\hat{\chi}_{E_t}(z) = \begin{cases} \frac{a_1 a_2}{2\pi\sqrt{(a_1 z \cdot v_1)^2 + (a_2 z \cdot v_2)^2}} J_1(\sqrt{(a_1 z \cdot v_1)^2 + (a_2 z \cdot v_2)^2}\sqrt{t}) & z \neq 0 \\ \frac{a_1 a_2 t}{4\pi} & z = 0 \end{cases}$$
$$(2.3)$$

**Lemma 3.** *We have*

$$A(R) = \sum_{n \neq 0} \frac{a_1 a_2}{\pi[(a_1 n \cdot v_1)^2 + (a_2 n \cdot v_2)^2]} J_2(\sqrt{(a_1 n \cdot v_1)^2 + (a_2 n \cdot v_2)^2}\sqrt{R}), \quad (2.4)$$

*the series converging uniformly and absolutely.*

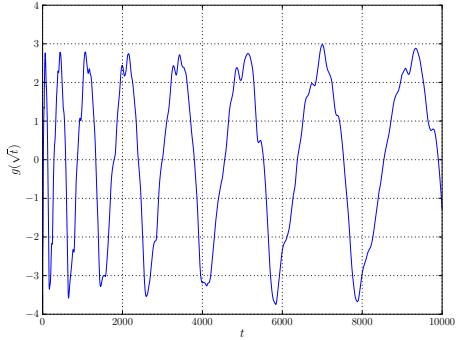 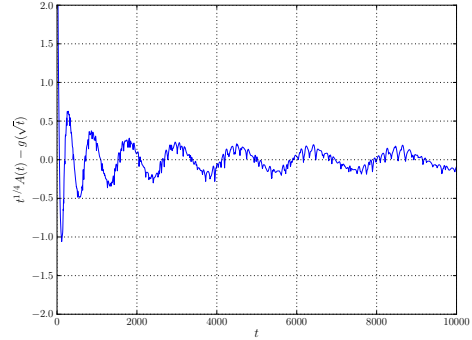

Figure 11: $g(\sqrt{t}), a_1 = 2, a_2 = 1/2$      Figure 12: $t^{1/4}A(t) - g(\sqrt{t})$

*Proof.* The same proof as for Lemma 1, with (2.3) used in place of (1.6).  □

**Theorem 4.** *The asymptotic expansions (1.15) and (1.16) hold, where now*

$$g(x) = -\frac{\sqrt{2}}{\pi^{3/2}} a_1 a_2 \sum_{n \neq 0} \left((a_1 n \cdot v_1)^2 + (a_2 n \cdot v_2)^2\right)^{-5/4}$$
$$\cdot \cos\left(\left((a_1 n \cdot v_1)^2 + (a_2 n \cdot v_2)^2\right)^{1/2} x - \frac{\pi}{4}\right). \quad (2.5)$$

*Proof.* Same as for Theorem 2, using Lemma 3 in place of Lemma 1.  □

See Figure 11 for $g(\sqrt{t})$ with $a_1 = 2$ and $a_2 = 1/2$ and Figure 12 for the difference $t^{1/4}A(t) - g(\sqrt{t})$ for the same family of ellipses. Plots for different families of ellipses are available on the website [5].

## 3  The Klein bottle and projective plane

If we identify the vertical boundaries of the square directly, and the horizontal boundaries with reflection, we obtain the standard flat Klein bottle KB. In terms of functions defined on the square, we are imposing the boundary conditions $u(0, y) = u(1, y)$ and $u(x, 0) = u(1 - x, 1)$ in order to have a function on KB. We may cover KB by the rectangular torus $[0, 1] \times [0, 2]$ with the identities

$$\begin{cases} u(x+1, y) = u(x, y) \\ u(1-x, y+1) = u(x, y) \end{cases} \quad (3.1)$$

describing the lifts of functions on KB to $\mathbb{R}^2$. The eigenfunctions of the Laplacian on KB lift to eigenfunctions on the rectangular torus, and so are linear combinations of functions of the form $e^{2\pi i\left(jx+\frac{k}{2}y\right)}$ with eigenvalue $(2\pi)^2(j^2 + (\frac{k}{2})^2)$. Now we observe that $e^{2\pi i(j(1-x)+\frac{k}{2}(y+1))} = (-1)^k e^{2\pi i(-jk+\frac{k}{2}y)}$. Thus there are two families of eigenfunctions

$$e^{2\pi i \frac{k}{2} y} \quad \text{for } k \text{ even (corresponding to } j = 0\text{), and} \tag{3.2}$$
$$e^{2\pi i\left(jx+\frac{k}{2}y\right)} + (-1)^k e^{2\pi i\left(-jx+\frac{k}{2}y\right)} \quad \text{for } j > 0. \tag{3.3}$$

We can therefore see that the eigenvalue function $N_{\text{KB}}$ is close to one half the counting function $N_{T_{1,2}}$ for the $[0,1] \times [0,2]$ torus.

**Theorem 5.** $N_{\text{KB}}(t) = \frac{1}{2} N_{T_{1,2}}(t) \pm \frac{1}{2}$.

*Proof.* $N_{T_{1,2}}(t)$ counts all integers $j, k$ such that $j^2 + (\frac{k}{2})^2 \leq \frac{t}{(2\pi)^2}$. When $j \neq 0$ the pair $\pm j$ contributes just a single eigenvalue to $N_{\text{KB}}(t)$. When $j = 0$ we count all $k$ such that $|k| \leq \frac{\sqrt{t}}{\pi}$ in $N_{T_{1,2}}(t)$, but just the even values of $k$ in $N_{\text{KB}}(t)$, and $\#\{k \text{ even} : |k| \leq \frac{\sqrt{t}}{\pi}\} = \frac{1}{2}\{k : |k| \leq \frac{\sqrt{t}}{\pi}\} \pm \frac{1}{2}$. $\square$

It is interesting to compare the Klein bottle with the projective plane (PP) obtained from $[0,1] \times [0,1]$ by identifying both sets of boundary edges with reflections. Functions on PP lift to $\mathbb{R}^2$ with the identities

$$\begin{cases} u(1-x, y+1) = u(x,y) \\ u(x+1, 1-y) = u(x,y) \end{cases} \tag{3.4}$$

and the torus $[0,2] \times [0,2]$ is a four-fold covering of PP. However, while it is possible to pull back the standard Laplacian to PP, the pairs $\{(0,0),(1,1)\}$ and $\{(0,1),(1,0)\}$ of identified points on PP are singularities (cone points with total angle $\pi$) with respect to the otherwise flat metric.

Reasoning as in the KB example, we know that eigenfunctions of the Laplacian on PP must be linear combinations of the functions $e^{2\pi i(\frac{j}{2}x+\frac{k}{2}y)}$ with eigenvalue $(2\pi)^2((\frac{j}{2})^2 + (\frac{k}{2})^2)$. Imposing the conditions (3.4) leads to four families of eigenfunctions:

$$\text{constants (corresponding to } j = 0 \text{ and } k = 0) \quad (3.5)$$

$$e^{2\pi i \frac{k}{2} y} + e^{-2\pi i \frac{k}{2} y} \text{ for } k > 0 \text{ even (corresponding to } j = 0 \text{ but } k \neq 0) \quad (3.6)$$

$$e^{2\pi i \frac{j}{2} x} + e^{-2\pi i \frac{j}{2} x} \text{ for } j > 0 \text{ even (corresponding to } k = 0 \text{ but } j \neq 0) \quad (3.7)$$

$$e^{2\pi i \left(\frac{j}{2} x + \frac{k}{2} y\right)} + e^{2\pi i \left(-\frac{j}{2} x - \frac{k}{2} y\right)} + (-1)^{j+k} \left( e^{2\pi i \left(-\frac{j}{2} x + \frac{k}{2} y\right)} + e^{2\pi i \left(\frac{j}{2} x - \frac{k}{2} y\right)} \right)$$
$$\text{for } j > 0 \text{ and } k > 0. \quad (3.8)$$

This leads to the identity

$$N_{\text{PP}}(t) = \frac{1}{4} N_{\text{T}_{2,2}}(t) + \frac{1}{4} \pm \frac{1}{2} = \frac{1}{4} N_{\text{T}_{1,1}}(4t) + \frac{1}{4} \pm \frac{1}{2}. \quad (3.9)$$

Similar results hold for KB and PP constructed from the tori considered in section 2. Related questions in the context of fractal Laplacians with the Sierpinski carpet replacing the square are discussed in [1].

563 MALOTT HALL, CORNELL UNIVERSITY, ITHACA, NY 14853, USA
*Email addresses:* dsj36@cornell.edu, str@math.cornell.edu